\documentclass[11pt]{article}
\usepackage{amsmath}
\usepackage{amssymb}
\usepackage{dsfont, mathrsfs}
\usepackage{amsthm}
\usepackage{arydshln}
\usepackage[colorlinks=true]{hyperref}
\usepackage{graphicx}
\DeclareGraphicsExtensions{.pdf,.png}
\graphicspath{{./figures/}{./}}
\usepackage{subfigure}

\usepackage{tikz}
\usetikzlibrary{calc}

\newcommand{\bR}{\mathbb R}
\def\a{\mathbf{a}}

\def\al{\boldsymbol{\alpha}}

\newtheorem{theorem}{\hspace*{\parindent}Theorem}

\newtheorem{conjecture}{Conjecture}

\newtheorem{question}[theorem]{Question}
\newtheorem{remark}{Remark}

\newcommand{\prth}[1]{\!\left( #1 \right) }

\newcommand{\Prob}[1]{ \mathbb{P} \prth{ #1 } } 

\def\longlongrightarrow{\hspace{+0.1ex} - \hspace{-1.1ex} - \hspace{-1.1ex} - \hspace{-1.1ex}\longrightarrow  } 

\newcommand{\tendvers}[2]{ \underset{#1 \rightarrow #2}{\longlongrightarrow} } 

\newcommand{\bracket}[1]{\left\langle #1 \right\rangle}

\def\Ai{{\operatorname{Ai}}}

\def\cP{\mathcal{P}}

\usepackage{authblk}

\begin{document}

\title{Open problems on polynomials, their zero distribution and related questions: 2023}
\author{Liudmyla\:Kryvonos$^{\rm }$\footnote{E-mail: L. Kryvonos -- \emph{liudmyla.kryvonos@vanderbilt.edu}}}

\affil{Vanderbilt University}
\date{}
\maketitle 
\thispagestyle{empty}	

\begin{abstract} 
	This paper collects open problems that were presented at the ``Hausdorff Geometry of Polynomials" workshop held on July 10-14, 2023 in Sofia, Bulgaria. 
\end{abstract}

\section*{Introduction}

 $\;\;\;\;\;$  Hausdorff Geometry of Polynomials workshop was held as part of the Mathematics Days in Sofia conference organized by the Institute of Mathematics and Informatics of the Bulgarian Academy of Science.\\

\textit{Organizers:} P. Boyvalenkov, A. Martinez-Finkelshtein,  E.B. Saff, and B. Shapiro.\\

\textit{Participants:} Y. Barhoumi-Andr\'eani, P. Boyvalenkov, D. Karp, V. Kostov, L. Kryvonos, A. Kuijlaars, A. Martinez-Finkelshtein, G. Nikolov, R. Orive, E.B. Saff, H. Sendov, B. Shapiro, N. Stylianopoulos, F. Wielonsky.

\vspace{0.6in}
In each section, we list an open problem, the corresponding proposer's name and their contact information.

\begin{center}
	\textbf{\Large{Yacine Barhoumi-Andr\'eani$^{\rm}$\footnote{E-mail: Y. Barhoumi-Andr\'eani -- \emph{yacine.barhoumi@gmail.com}}}}\\
	
\end{center}

\begin{center}
	\textbf{\large{Some open problems in Random Matrix Theory and Orthogonal Polynomials}}
\end{center}
\vspace{0.1in}

\paragraph{\large{1. A problem on the  GUE Tracy-Widom distribution}}\mbox{}\\

The GUE Tracy-Widom distribution is the probability distribution defined by the following probability distribution function~:
\begin{align*}
	\Prob{ \mathcal{TW}_2 \leq s} := \det\prth{ I - \boldsymbol{K}_{\! \Ai} }_{ L^2( [s, +\infty) ) }
\end{align*}
where $ \boldsymbol{K}_{\! \Ai} $ is the Airy kernel operator acting on functions $ f \in L^2([s, +\infty)) $ by $ \boldsymbol{K}_{\! \Ai} f(x) := \int_{[s, +\infty)} K_{\! \Ai}(x, y) f(y) dy $ and
\begin{align*}
	K_{\! \Ai}(x, y) := \int_{\mathbb{R}_+}\Ai(x + t) \Ai(y + t) dt
\end{align*}

Here, $ \Ai $ is the Airy function, i.e. the bounded solution of the ODE $ y'' = xy $. The Airy kernel operator can be understood as the $ L^2(\mathbb{R}_+) $-square of the Hankel operator $ \boldsymbol{H}(\Ai) $ whose kernel is $ (x, y) \mapsto \Ai(x + y) $. Writing 
\begin{align*}
	K_{\! \Ai}(x, y) = \int_s^{+\infty} \Ai(x + t - s) \Ai(y + t - s) dt
\end{align*}
one sees that $ \boldsymbol{K}(\Ai)_{\curvearrowright L^2( [s, +\infty) ) } = \prth{ \boldsymbol{H}(\Ai(\cdot - s))_{\curvearrowright L^2( [s, +\infty) ) } }^2 $ i.e. one also has a square on $ L^2( [s, +\infty) ) $.

\medskip

In their foundational article that now defines the Tracy-Widom distribution \cite{TracyWidomCMP}, the authors define in \cite[\S~IV-B p. 165]{TracyWidomCMP} a self-adjoint differential operator on $ L^2( [s, +\infty) )$ that commutes with $ \boldsymbol{H}(\Ai(\cdot - s)) $ (hence to its square $ \boldsymbol{K}_{\!\Ai} $). This differential operator is given by 
\begin{align*}
	\mathcal{L}_{TW, s} := \frac{d}{dx} \alpha_s  \frac{d}{dx} + \beta_s ,   \qquad \alpha_s(x) := x - s, \quad \beta_s(x) := -x(x - s)
\end{align*}

This situation mimics the case of the ``bulk'' in Random Matrix Theory where the equivalent of $ \boldsymbol{K}_{\!\Ai} $ is given by the \textit{sine kernel operator} $ \boldsymbol{K}_{\!\operatorname{sinc}} : f \mapsto f * \operatorname{sinc} $ where the convolution is taken on $ L^2([-1, 1]) $ or $ L^2([-c, c]) $ ($ c = \pi $ and sinc is replaced with $ \operatorname{sinc}(c\cdot) $ then). The commuting differential operator is then given by a famous result of Slepian by the \textit{prolate spheroidal wave function} operator $ \operatorname{PSO}_c := \frac{d}{dx} (c^2 - x^2) \frac{d}{dx} - c^2 x^2 $ whose eigenvectors are well studied and can be expressed as linear combinations of Legendre functions.

\medskip

Such a study was, to the best of my knowledge, never performed for $ \mathcal{L}_{TW, s} $, despite its ``simple'' polynomial structure (it is not of the Heun or hypergeometric type, and the fact that $ \beta_s $ is a quadratic polynomial prevents a priori from using a polynomial expansion).

\medskip

It is of interest due to the Lidskii formula that states that for a trace-class operator $ \boldsymbol{K} $ acting on $ L^2(A) $ for a certain set $A$, one has
\begin{align*}
	\det(I - \boldsymbol{K})_{L^2(A)} = \prod_{k \geqslant 1} (1 - \lambda_k(A))
\end{align*}
where the $ \lambda_k(A) $ are the eigenvalues of $ \boldsymbol{K} \curvearrowright L^2(A) $. 

The commuting differential operator technique allows to get the eigenvectors $ \psi_k $, hence the eigenvalues $ \bracket{\boldsymbol{K} \psi_k, \psi_k}_{\! L^2(A) } $ (if $ |\!|\psi_k|\!|^2 = 1 $). As a result, one could write the Tracy-Widom distribution function as a Lidskii product with several probabilistic applications (for instance, Tracy and Widom used the WKB method to have access to the tail behaviour of their distribution, etc.).

Overall, since the universality class of the GUE Tracy-Widom distribution (the so-called ``KPZ universality class'') is now very huge (lots of models have been proven to converge to this distribution), these eigenvectors must have several interesting properties. The Airy kernel operator would write (with obvious notations)
\begin{align*}
	K_{\Ai, s} = \sum_{k \geqslant 0} \lambda_k(\Ai, s) \psi_{k, s} \otimes \psi_{k, s}, 
	\qquad 
	\lambda_k(\Ai, s) & = \left\langle \boldsymbol{K}_{\!\Ai} \psi_{k, s}, \psi_{k, s} \right\rangle_{\! L^2([s, +\infty))} \\
	& = \vert\!\vert \boldsymbol{H}(\Ai)\psi_{k, s}\vert\!\vert_{L^2([s, +\infty) )}^2
\end{align*}

To sumarise~:
\begin{question}
	Study the eigenvectors of the operator
	\begin{align*}
		\mathcal{L}_{TW, s} := \frac{d}{dx} \alpha_s  \frac{d}{dx} + \beta_s ,   \qquad \alpha_s(x) := x - s, \quad \beta_s(x) := -x(x - s)
	\end{align*}
\end{question}

\medskip

\medskip
\paragraph{\large{2. A related problem}}\mbox{}\\

The GUE kernel operator is the restriction to $ L^2([s, +\infty)) $ of the projection on the first $n$ vectors of $ L^2(\mathbb{R})$ with the (probabilistic) Hermite functions, i.e.
\begin{align*}
	K_n(x, y) := \sum_{k = 0}^{n - 1} \phi_k(x)\phi_k(y), \qquad \phi_k(x) := \frac{H_k(x)}{\sqrt{k!}} \frac{e^{-x^2/4} }{(2\pi)^{1/4}} \quad H_k = \mbox{Hermite polynomials}
\end{align*}

This is the rescaling of this operator ``at the edge'' that gives the Airy kernel operator, i.e.
\begin{align*}
	\sigma_n K_n(\mu_n + \sigma_n x, \mu_n + \sigma_n y) \tendvers{n}{+\infty} K_{\!\Ai}(x, y), \quad  \sigma_n := n^{-1/6}, \qquad \mu_n := 2\sqrt{n}
\end{align*}

One question would thus be~:
\begin{question}
	Can one find a commuting differential operator to $ \boldsymbol{K}_{\! n} $ acting on $ L^2([s, +\infty)) $~? What are its eigenvectors~?
\end{question}

Such an operator is not necessarily of the second order, but it should converge in some sense to $ \mathcal{L}_{TW, s} $ when suitably rescaled.


\bibliographystyle{amsplain}

\vspace{0.8in}

\begin{center}	
	\textbf{\Large{Dmitrii Karp$^{\rm}$\footnote{E-mail: D. Karp -- \emph{dimkrp@gmail.com}}}}
\end{center}	

\paragraph{\large{1.Toeplitz determinants of Pochhammers.}}\mbox{}\\

For a given finite numerical sequence $\{f_{k}\}_{k=0}^{n}$ define the polynomial $P_n(z)$ by
$$
P_n(z)\!=\!\sum\limits_{k=0}^{n}f_{k}f_{n-k}\binom{n}{k}\left[(z)_k(z)_{n-k}-(z+1)_k(z-1)_{n-k}\right],
$$
where $(z)_0=1$, $(z)_{k}=z(z+1)\cdots(z+k-1)$, $k\ge1$, is the Pochhammer symbol.

\bigskip
\textbf{Conjecture~1.} \textit{Suppose the polynomial $\sum\nolimits_{k=0}^{n}f_kz^k$ has only real negative zeros. Then the same is true for the polynomial $P_n(z)$ defined above.}

\bigskip

Known: the polynomial $P_n(z)$ has degree $n-2$ and  positive coefficients  \cite[Theorem~1]{K2013}.

\bigskip

A similar conjecture can be also be stated for polynomials formed by higher order Toeplitz determinants, namely,
\begin{multline*}
	P_n^r(z):=\!\!\!\!\!\!\!\!\sum\limits_{k_1+k_2+\cdots+k_r=n}
	\!\!\binom{n}{k_1,k_2,\ldots,k_r}f_{k_1}f_{k_2}\cdots{f_{k_r}}
	\\
	\times\begin{vmatrix}
		(z)_{k_1} & (z+1)_{k_1} &\cdots & (z+r-1)_{k_1}\\
		(z-1)_{k_2} & (z)_{k_2} &\cdots & (z+r-2)_{k_2}\\
		\vdots & \vdots &\cdots &\vdots\\
		(z-r+1)_{k_r} & (z-r+2)_{k_r} &\cdots & (z)_{k_r}\\
	\end{vmatrix}.
\end{multline*}

It is probably not too hard to prove that the polynomial $P_{n}^{r}(z)$ has degree $n-r(r-1)$. 

\medskip

\textbf{Conjecture~2.}
\textit{Suppose the polynomial $\sum\nolimits_{k=0}^{n}f_kz^k$ has only real negative zeros. Then the same holds for the polynomial $P_{n}^{r}(z)$ for each $r\geq{2}$.}

\bigskip

One more conjecture about $P_{n}^{r}(z)$ which seem to hold numerically is the following.

\smallskip

\textbf{Conjecture~3.}
\textit{Suppose $\{f_k\}_{k=0}^{n}$ is a finite \textit{P\'{o}lya frequency sequence} of order  $d\geq{r}$. Then the
	polynomial $P_{n}^{r}(z)$ is Hurwitz
	stable (all its zeros lie on the open left half-plane). In particular, it has positive coefficients.}

Further details can be found in \cite{K2013}.

\paragraph{\large{2. Hypergeometric generating function.}}\mbox{}\\

Suppose $a_1,a_2,\ldots,a_d,b>0$ are real and $m_1,m_2,\ldots,m_d>0$ are integer. Define the polynomial $A(x)$ by 
\begin{equation*}
	{_{d+1}F_{d}}\!\left(\begin{array}{c}a_1+m_1,a_2+m_2,\ldots,a_d+m_d, b\\a_1,a_2,\ldots,a_d\end{array}\vline\:x\right)=\frac{A(x)}{(1-x)^{m+b}},
\end{equation*}
where $m=\sum_{j=1}^{d}{m_j}$.  One can see that $A(x)$ is indeed a polynomial via the so-called Miller-Paris transformation \cite[Theorem~2]{KP2020}.

\bigskip

\textbf{Open problem.}  \textit{Find sufficient conditions on $a_1,a_2,\ldots,a_d,b,m_1,m_2,\ldots,m_d$ ensuring that $A(x)$ has only  real zeros (only negative zeros). }

\bigskip

\textit{Example}: $d$-Narayana polynomial can be defined by
$$
{_{d}F_{d-1}}\left(\left.\!\!\begin{array}{c}m+1,m+2,\ldots,m+d\\2,3,\ldots,d\end{array}\right|x\!\right)=\frac{N_{d,m}(x)}{(1-x)^{dm+1}}
$$
and is known to have only real and negative zeros.

\medskip

Related problems for the entire hypergeometric functions ${}_pF_{d}$ have been better studied.  In particular, given that $a_1,a_2,\ldots,a_d>0$ the function
$$
{_{d}F_{d}}\!\left(\begin{array}{c}a_1+m_1,a_2+m_2,\ldots,a_d+m_d\\a_1,a_2,\ldots,a_d\end{array}\vline\:x\right)
$$
has only real (and hence negative) zeros  if and only if $m_1,m_2,\ldots,m_d$ are positive integers. While the function ($p<d$)
$$
{_{p}F_{d}}\!\left(\begin{array}{c}a_1+m_1,a_2+m_2,\ldots,a_{p}+m_p\\a_1,a_2,\ldots,a_d\end{array}\vline\:x\right)
$$
also has only real negative zeros if $m_1,m_2,\ldots,m_d$ are positive integers, but this condition is not necessary. Finding other sufficient conditions remains an open problem.

\vspace{0.8in}
\begin{center}
	\textbf{\Large{Arno Kuijlaars$^{\rm}$\footnote{E-mail: A. Kuijlaars -- \emph{arno.kuijlaars@kuleuven.be}}}}\\
	
\end{center}

\begin{center}
	\textbf{\large{Orthogonal polynomials in the complex plane with varying
	exponential weights and their limiting zero counting measures}}
\end{center}

Let $P$ be a monic polynomial of degree $n$ with zeros
$z_{1}, \ldots, z_{n}$, listed according multiplicities.
We use
\[ \nu(P) = \frac{1}{n} \sum_{j=1}^n \delta_{z_j} \]
to denote the normalized zero counting measure of $P$. 
An external field on $\mathbb C$ is a continuous function 
$V : \mathbb C \to \mathbb R \cup \{ +\infty \}$, not identically $+\infty$, with 
\begin{equation} \label{Vgrowth}
	\lim_{|z| \to \infty} \left( V(z) - \log (1+|z|) \right) = + \infty.
\end{equation}
The growth condition \eqref{Vgrowth} guarantees that, for each $n$,
there is a unique monic polynomial $P_n$ of degree $n$ such that
\begin{equation} \label{Pnorthog}
	\int_{\mathbb C}
	P_n(z) \overline{z^k} e^{-n V(z)} dA(z) =
	0, \qquad k =0,1, \ldots, n-1 \end{equation}
where $dA$ denotes the area measure in the complex plane. \\

\textbf{Problem.}
\textit{Find  external fields  $V$ such that the sequence of normalized
	zero counting measures of the orthogonal polynomials \eqref{Pnorthog} 
	has a weak limit, say $\nu(P_n) \stackrel{\ast}{\longrightarrow} \mu_V$ as $n \to \infty$,
	and characterize the limiting measure $\mu_V$.}\\

The problem is trivial for radial external fields $V$, 
i.e., if $V(z) = V(|z|)$ for every $z \in \mathbb C$, since in such a
case $P_{n}(z) = z^n$ for every $n \in \mathbb N$, and $\mu_V$ is the
Dirac point mass at $z=0$.

The analogous problem on the real has been very well studied. In that
case the limiting measure $\mu_V$ exists and it is the equilibrium
measure in the external field $V$. That is, $\mu_V$ is the unique
probability measure that minimizes
\[ \iint \log \frac{1}{|x-y|} d\mu(x) d\mu(y)
+ \int V(x) d\mu(x) \]
among all probability measures $\mu$ on the real line,
see e.g.\ \cite{SaTo}.

This result is not true in the complex plane. Indeed, the 
support of the equilibrium measure in an external field on $\mathbb C$ 
is  a planar domain, see \cite{ElFe, SaTo},  such as a disk or an annulus 
in case $V$ is radial.
Thus for radial external fields the limit
of the normalized counting measures does not coincide with 
the equilibrium measure in the external field.

External fields of the form $V(z) = |z|^2 -  2c \log |z-a|$ 
with $c > 0$ and $a \in \mathbb C \setminus \{0\}$ were studied
in \cite{BBLM} and it was found that the limit $\mu_V$ of the normalized
zero counting measures exists and $\mu_V$ is supported on a contour that
lies inside the convex hull of the support of the equilibrium measure.
It is expected that this is a general phenomonon. It would be very
interesting to identify larger classes of external fields for
which this hapens.

\vspace{0.8in}
\begin{center}
	\textbf{\Large{Geno Nikolov$^{\rm }$\footnote{E-mail: G. Nikolov -- \emph{geno@fmi.uni-sofia.bg}}}}\\
	
\end{center}

\begin{center}
	\textbf{\large{Snake polynomials with positive Chebyshev
			expansions}}
\end{center}

Denote by $\cP_m$ the set of all algebraic polynomials of degree at most
$m$. Without loss of generality, we can assume that coefficients of the polynomials are real.

\begin{figure}[htp]
	\begin{center}
		\centering
		\subfigure[]{\includegraphics[width=0.33\textwidth]{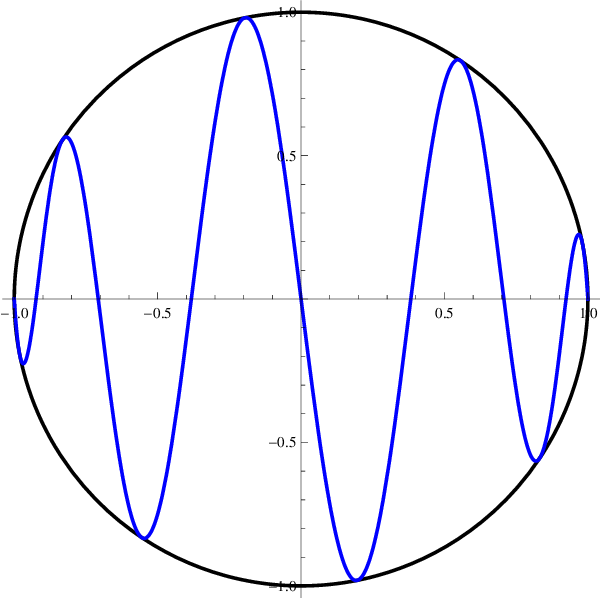}} $\;\;\;$
		\subfigure[]{\includegraphics[width=0.33\textwidth]{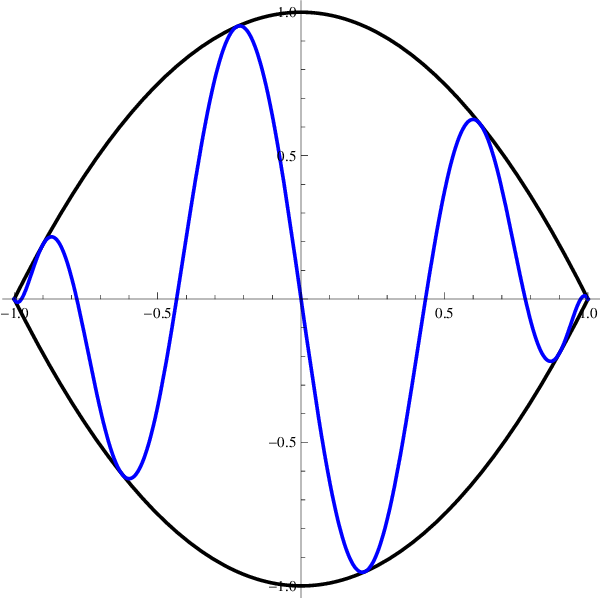}} 
		\subfigure[]{\includegraphics[width=0.33\textwidth]{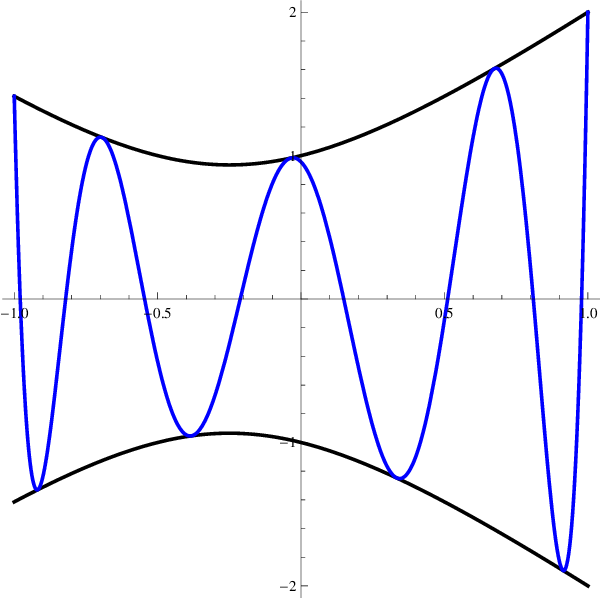}}$\;\;\;\;\;$
		\subfigure[]{\includegraphics[width=0.33\textwidth]{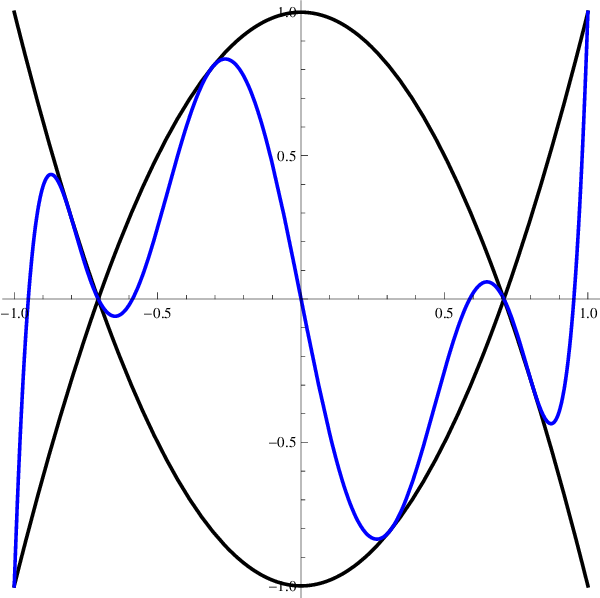}}
		\caption{The snake polynomials for various majorants: (a) $\mu(x)=\sqrt{1-x^2} \;\;$  (b) $\mu(x)=1-x^2\;\;$  (c) $\mu(x)=\sqrt{2x^2+x+1}\;\;$  (d) $\mu(x)=|1-2x^2|$}
		\label{fig:Snake}
	\end{center}
\end{figure}

Consider a majorant $\mu\in C[-1,1],\ \mu(x)\geq 0, \ x\in [-1,1]$. If there exists $P\in\cP_n$, $P\ne 0$, such that $-\mu(x)\leq
P(x)\leq\mu(x)$, $x\in [-1,1]$, then there exists a unique (up to
orientation) polynomial $\omega_{\mu}\in \cP_n$ called \textit{snake polynomial associated with} $\mu$,
which oscillates most between $\pm \mu$ (see Fig \ref{fig:Snake}).   

The $n$-th snake polynomial
$\omega_{\mu}$ is uniquely determined by the following properties:
\medskip

a) $\quad |\omega_{\mu}(x)|\leq\mu(x)\,\quad x\in [-1,1]\,$;
\medskip

b)~~ There exists a set $\delta^{*}=(\tau_i^{*})_{i=0}^{n}$,
$1\geq\tau_0^{*}>\cdots>\tau_n^{*}\geq -1$, such that
$$
\omega_{\mu}(\tau_i^{*})=(-1)^{i}\mu(\tau_i^{*}),\quad i=0,\ldots,n,
$$

where $\delta^{*}$ is referred to as  \textit{the set of alternation points
	of} $\omega_{\mu}$.

\newpage

\textbf{Open problem.} \textit{Find a class of majorants $\mu$ for which
the associated snake polynomials $Q_n$ have non-negative or sign alternating
expansion in the Chebyshev polynomials of the first kind, i.e.
$$
Q_n=\sum_{i=0}^{n}a_i\,T_i\ \text{ with all }\ a_i\geq 0\ \text{ or
	with }\ a_i\,a_{i+1}<0\ \ \forall\ i.
$$}

\textbf{Conjecture (G.N.)} \textit{If $\mu\geq 0$ is a continuous even convex function in $[-1,1]$,
then the associated with $\mu$ snake polynomials have non-negative
expansion in the Chebyshev polynomials of the first kind.}

The snake polynomials are analogous to the Chebyshev polynomials $T_{n}$ for $\mu\equiv 1$ and are of great interest in relation to the following two problems:

\textbf{Problem~1: Markov inequality with a majorant}
Given $n,k \in \mathbb{N}$, $1\le k\le n$, and a majorant $\mu \ge
0$, find
$$
M_{k,\mu} := \sup\{\|p^{(k)}\|\,:\, p\in\cP_n,\,|p(x)| \le \mu(x),\,
x\in[-1,1]\}\, .
$$

\textbf{Problem 2: Duffin-Schaeffer inequality with a majorant}
Given $n,k \in \mathbb{N}$, $1\le k\le n$, and a majorant $\mu \ge
0$, find
$$
D_{k,\mu}^{*} := \sup\{\|p^{(k)}\|\,:\, p\in\cP_n,\,|p(x)| \le
\mu(x),\, x\in\delta^{*}\}\, .
$$

A conjecture (belonging to mathematical folklore) says that the
extremal polynomial to Problem 1 is the snake polynomial
$\omega_{\mu}$. So far, no counterexample to this conjecture has been
found. On the contrary, $\omega_{\mu}$ is not always the extremal
polynomial to Problem 2, the following counterexamples are known:
$$
1)\quad\mu(x)=\sqrt{1-x^2},\quad k=1;\qquad
2)\quad\mu(x)=1-x^2,\quad k=1,2\,.
$$

In \cite{NS1, NS} the majorants $\mu$ were studied for which the snake-polynomial $\omega_{\mu}$
	is extremal to both Problems 1 and 2, i.e., when equality
	\begin{equation}\label{eq}
		M_{k,\mu}   =
		D_{k,\mu}^* = 
		\|\omega^{(k)}_\mu\|\,  
	\end{equation}
holds. It was proved that (\ref{eq}) holds if $\omega_{\mu}$ admits non-negative or sign-alternating expansion in the Chebyshev polynomials, and so the question of our interest is to characterize the corresponding class of measures $\mu$.

\vspace{0.8in}
\begin{center}
	\textbf{\Large{Edward Saff$^{\rm}$\footnote{E-mail: E. Saff -- \emph{edward.b.saff@vanderbilt.edu}}}}\\
	
\end{center}

\begin{center}
	\textbf{\large{Equilibrium measure on an annulus}}
\end{center} 

For $-2<s< d$, the Riesz energy of a probability measure $\mu$ on $\mathbb{R}^{d}$ is
$$
I_{s} (\mu):= \iint K_{s}(x-y) d\mu(x) d\mu(y), 
$$
 where $K_{s}(x)$ is the Riesz $s$-kernel defined by $K_{s}(x) = \begin{cases}
 	\frac{1}{s|x|^{s}} ,\;\;\; \;\;\; \;\;\textnormal{if} \;\;s \ne 0,\\
 	 - \log|x|, \;\;\; \textnormal{if}\;\; s = 0.
 	\end{cases}
 $
  
  \textbf{Riesz original problem.} In \cite{Riesz} Riesz showed that the equilibrium measure 
  
   $\mu_{*} :=   \mathop{\arg \;\; \min_{\mu}}_{\textnormal{supp}(\mu) \subset B_{R}} \;\;  I_{s}(\mu)$ on the ball $B_{R}:=\{x \in \mathbb{R}^{d}: |x| \leq R\}$ is given by 
  $$ \mu_{*} =
  \begin{cases}
  	\sigma_{R}, \hspace{1.5in} \:\quad \quad \quad  \textnormal{if}  -2<s\leq d-2,\\
  	\frac{\Gamma(1+\frac{s}{2})}{R^{s} \pi^{\frac{d}{2}} \Gamma(1+\frac{s-d}{2})} \frac{\textbf{1}_{B_{R}}}{(R^{2} - |x|^{2})^{\frac{d-s}{2}}} dx, \; \;\;\;\; \quad  \textnormal{if}  \;\; 0\leq  d-2 <s < d,\\
  	\end{cases}
  $$
  where $\sigma_{R}$ is the uniform distribution on sphere of radius $R$. 
  
  The math arXiv version (arXiv:2108.00534) of \cite{CSW} contains in its Appendix D.2 a short alternative proof of the Riesz original problem.
  
 \textbf{Equilibrium open problem for multi-dimensional annulus.} \textit{Given $d-2 < s < d$}, find an equilibrium measure on an annulus $B_{r,R}:=\{x \in \mathbb{R}^{d} : r\leq|x| \leq R\}$, i.e. the measure $\mu_{*}$ such that
 $$
 \mu_{*} =   \mathop{\arg \;\; \min_{\mu}}_{\textnormal{supp}(\mu) \subset B_{r,R}} \;\;  I_{s}(\mu).  
 $$

\vspace{1in}
\begin{center}
	\textbf{\Large{Boris Shapiro$^{\rm}$\footnote{E-mail: B. Shapiro -- \emph{shapiro@math.su.se}}}}\\
	
\end{center}

\paragraph{\large{1. On shadows of polynomials}}\mbox{}\\

In \cite{BoHaSh} my coauthors and I studied asymptotic of root distribution for the polynomial sequences $\{Q_{n,\alpha}(x)\}_{n=1}^\infty$ where 
$$Q_{n,\alpha}(x)=\frac {d^{[\alpha n]}P^n(x)}{dx^{[\alpha n]}}. 
$$
Here $P(x)$ is a fixed univariate polynomial of degree $d\ge 2$, $\alpha \in (0,d)$ and $[\alpha n]$ stands for the integer part of $\alpha n$. 

\begin{figure}[htp]
	\begin{center}
		\includegraphics[width=.48\textwidth]{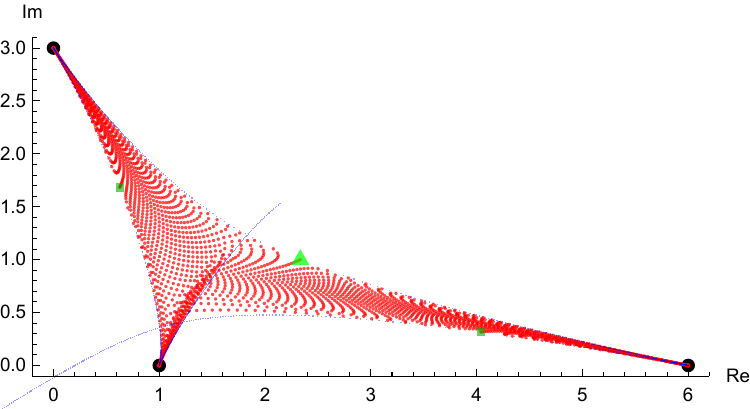}\hfill
		\includegraphics[width=.48\textwidth]{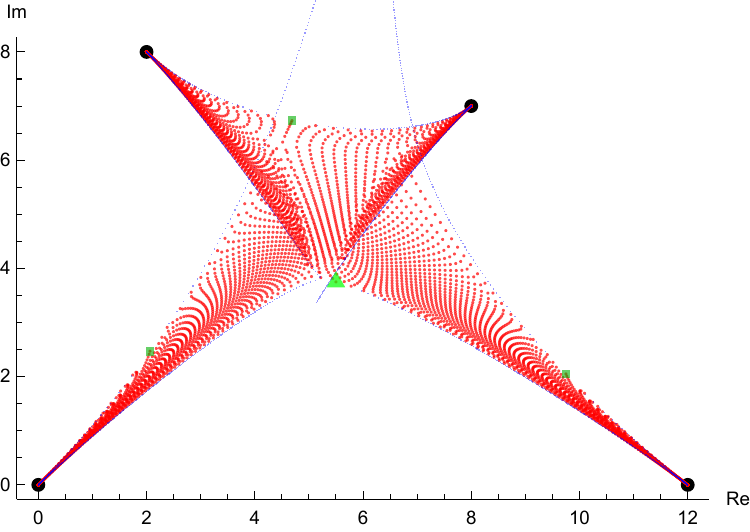}
		\includegraphics[width=.48\textwidth]{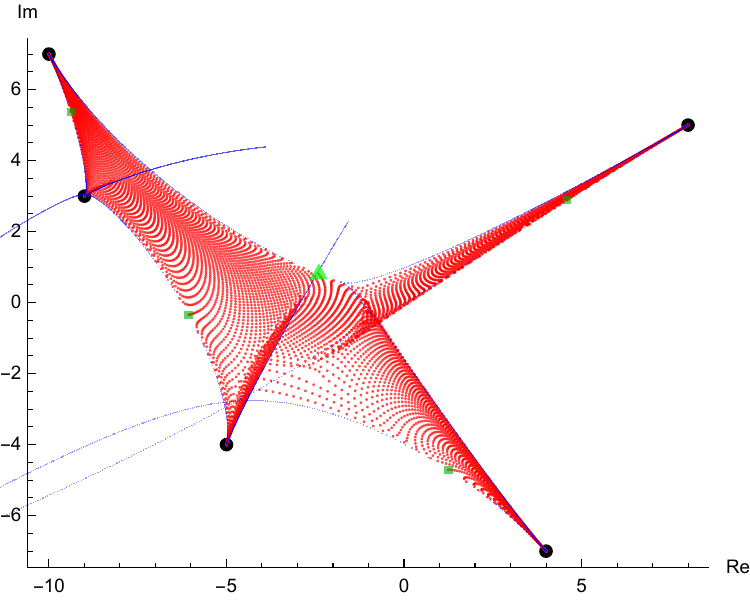}\hfill
		\includegraphics[width=.48\textwidth]{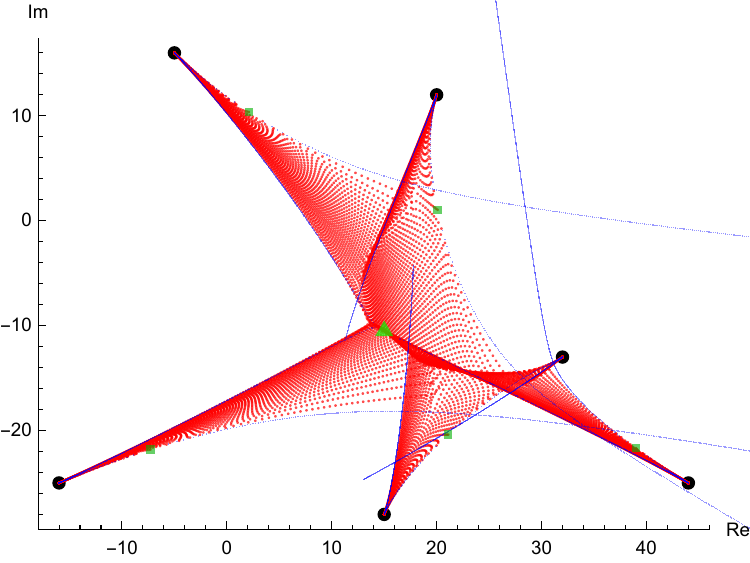}
	\end{center}
	\caption{The union of all zeros of $Q_{30,m}(z)$ for $m= 0,1,$ $\dotsc,30\deg{P}-1$ shown by small red dots. The large black dots are the zeros of $P$, the green squares are the critical points of $P$, the green triangle is the center of mass of the zero locus of $P$.}
	\label{fig:polShadowExamples}
\end{figure}

Motivated by the set-up of \cite{BoHaSh} and the famous Gauss-Lucas theorem, let us consider the double-indexed polynomial family
$$Q_{n,m}(x)=\frac {d^{m}P^n(x)}{dx^{m}}, 
$$ 
where $0\le m \le nd$. Obviously, for any pair of non-negative integers  $(n,m)$, all roots of $Q_{n,m}(x)$ (if any) lie in the convex hull of the roots of $P(x)$. Consider the union of all roots of all  $Q_{n,m}(x)$ with $0\le m \le nd$ and denote the closure of the latter countable set as $\Upsilon_P$ -- the {\it shadow of the polynomial} $P(z)$. Numerical examples of  shadows for polynomials of different degrees are shown in Fig.~\ref{fig:polShadowExamples}. 

\begin{conjecture}\label{conj:main}For any polynomial $P$ of degree at least $3$ whose roots are in convex position, but do not form a regular polygon,
	
	\noindent
	{\rm(i)}  $\Upsilon_P$ is a closed domain in the convex hull of the roots of $P$.
	
	\noindent
	{\rm(ii)} All critical points of $P$ lie on the boundary of $\Upsilon_P$.
	
	\noindent
	{\rm(iii)} the boundary of $\Upsilon_P$ has no inflection points.
	
	\noindent
	{\rm(iv)} The boundary of $\Upsilon_P$ is contained in the union of all critical values (w.r.t. $z$) of the rational function
	\begin{equation}
		F_{\al(z)}=z-\al \frac{P(z)}{P^\prime(z)}.
	\end{equation}
	where the parameter $\al$ runs over the interval $[0,d]$.  
\end{conjecture}

\begin{remark} Item (iii) can be reformulates as follows. The boundary of  $\Upsilon_P$ consists of all $u$ for which the family $\tilde \Phi(\al,z,u)=\al P(z)+(u-z)P^\prime(z)$ has a multiple root w.r.t. $z$ when $\al\in [0,d]$.
\end{remark} 

It seems that the shadow  $\Upsilon_P$  is a very natural and interesting domain on whose boundary lie the critical points of $P$. It can also considered as the support of an  asymptotic root-counting measure which is worth studying in its own rights.  

\paragraph{\large{2. On Maxwell's problem}}\mbox{}\\ 

In \cite{GNSh} my coauthors and I have been discussing  the following ''result" which they found in \cite {Ma}. 

\begin{conjecture}[J.~C.~Maxwell]\label{conj:Max} For any configuration of $N$ fixed  electric charges in $\bR^3$ having finitely many points of equilibrium, the number of the latter points is at most $(N-1)^2$.
\end{conjecture} 

\begin{remark} Observe that almost all configurations of fixed point charges satisfy the assumptions of the above conjecture that their number of points of equilibrium is finite. 
\end{remark} 

The following folklore claim is widely open. 

\begin{conjecture} Any configuration of $N$ positive  electric charges in $\bR^d$ has only finitely many points of equilibrium. 
\end{conjecture} 

In \cite{GNSh} we were only able to prove  a very bad  upper bound for the number of points of equilibrium however in a more general situation. Since then our upper bound has been improved in a number of very special cases. The best general improvement has been made in a recent preprint \cite{Zo}.  However Conjecture~\ref{conj:Max} is still open even for $n=3$.

\end{document}